\documentclass{amsart}
\usepackage{amsmath,amssymb,amsthm,amscd}
\usepackage{pgfplots}
\usepackage{mathrsfs}
\newtheorem{formula}{}[section]

\newtheorem{lemma}[formula]{Lemma}
\newtheorem{theorem}[formula]{Theorem}
\theoremstyle{definition}
\newtheorem{definition}[formula]{Definition}
\newtheorem{example}[formula]{Example}
\theoremstyle{remark}
\newtheorem*{remark}{Remark}

\begin{document}

\title{Narrow positively graded Lie algebras}
\author{Dmitry V. Millionshchikov}
\thanks{This work is supported by the Russian Science Foundation under grant 14-11-00414.}
\address{Steklov Institute of Mathematics of RAS, 8 Gubkina St. Moscow,
119991, Russia and  Department of Mechanics and Mathematics, Moscow
State University, 1 Leninskie gory, 119992 Moscow, Russia}
\email{million@mech.math.msu.su}

\begin{abstract}
We classify real and complex infinite-dimensional positively graded Lie algebras
${\mathfrak g}=\oplus_{i=1}^{{+}\infty}{\mathfrak g}_i$  with properties
$$
[{\mathfrak g}_1, {\mathfrak g}_i]={\mathfrak g}_{i{+}1},  \; \dim{{\mathfrak g}_i}+\dim{{\mathfrak g}_{i{+}1}} \le 3, \; i \ge 1.
$$
In the proof of the main theorem we apply successive central extensions of finite-dimensional Carnot algebras.
In sub-Riemannian geometry, control theory, and geometric group theory, Carnot algebras play a significant role.
\end{abstract}
\date{}

\maketitle

\section*{Introduction}
Zelmanov and Shalev introduced in \cite{ShZ}  the concept of width for positively graded Lie algebras. A positively graded Lie algebra ${\mathfrak g}=\oplus_{ i \in {\mathbb N}} {\mathfrak g}_i$ is called a Lie algebra of width $d$, if there is such a number  $d$ (minimal with this property), that $\dim{{\mathfrak g}_i} \le d, \forall i \in {\mathbb N}$. The problem of classifying graded Lie algebras of finite width was defined by Zel'manov and Shalev as an important and difficult \cite{ShZ}.

The present paper is devoted to the classifications of narrow ("width 3/2") infinite-dimensional naturally graded Lie algebras. They can also be considered as Lie algebras of slow linear growth. In the theory of nonlinear hyperbolic partial differential equations, the concept of the characteristic Lie algebra of the equation is defined. Two graded Lie algebras from our list ${\mathfrak n}_1$ and ${\mathfrak n}_2$, which are positive parts of the affine Kac-Moody algebras $A_1^{(1)}$ and $A_2^{(2)}$ respectively are isomorphic to the characteristic Lie algebras of the sine-Gordon and Tzitzeika equations \cite{Mill3}. 

In conclusion, it should also be noted that problems related to narrow and slowly growing Lie algebras have been extensively studied in the case of a field of positive characteristic \cite{CMN_97}.

\section{Positively graded Lie algebras}

\begin{definition}
A Lie algebra ${\mathfrak g}$ is called ${\mathbb N}$-graded, if there exists a decomposition of it into a direct sum of linear subspaces ${\mathfrak g}_i$ such that
$$
{\mathfrak g}{=}\bigoplus_{i \in{\mathbb N}}{\mathfrak g}_i, \; [{\mathfrak g}_i, {\mathfrak g}_j] \subset  {\mathfrak g}_{i{+}j}, \;  i,j \in  {\mathbb N}.
$$
\end{definition}

\begin{definition}
It is said that a ${\mathbb N}$-graded Lie algebra ${\mathfrak g}=\oplus_{i \in{\mathbb N}}{\mathfrak g}_i$ has a width $d$, if there exists a positive number $d$ (minimal with this property), that 
$$
\dim{{\mathfrak g}_i} \le d, \;   \forall i \in  {\mathbb N}.
$$
\end{definition}

We give five examples of infinite-dimensional ${\mathbb N}$-graded Lie algebras of width one. For all these Lie algebras, the homogeneous subspaces ${\mathfrak g}_i$ are one-dimensional and in each of them a basis vector  $e_i, \; i \in {\mathbb N}$ is fixed. 

\begin{equation}
\label{5primerov}
\begin{array}{l}
1) \;\; \mathfrak{m}_0: \;\; \;
[e_1,e_i]=e_{i+1}, \; \forall \; i \ge 2;\\

2) \;\;\mathfrak{m}_2: \;\; 
\left\{ \begin{array}{c}
[e_1, e_i ]=e_{i+1},  \;\; i \ge 2, \\

[e_2, e_j ]=e_{j+2},  \; j \ge 3.
\end{array}
\right.\\

3) \;\; W^+:\; [e_i,e_j]= (j-i)e_{i{+}j}, \; \forall \;i,j \in {\mathbb N},\\

4) \;\; {\mathfrak n}_1: \;[e_i,e_j]= c_{i,j}e_{i{+}j}, \; c_{i,j}=\left\{ \begin{array}{c} 1, \; j{-}i \equiv 1 \; {\rm mod} \; 3;\\
0,  \; j{-}i \equiv 0  \; {\rm mod}\; 3;\\
{-}1,  \; j{-}i \equiv{-}1 \; {\rm mod}\; 3;\\
\end{array}\right.\\

5) \;\;  {\mathfrak n}_2: \; [e_q,e_l]= d_{q,l}e_{q{+}l}, \; q,l \in {\mathbb N},
\end{array}
\end{equation}
where the structure constants $d_{q,l}$ of the Lie algebra ${\mathfrak n}_2$ are defined by the Table  \ref{structure_const_n_2}.  
The matrix $(d_{q,l})$ is skew-symmetric, and its elements $d_{q,l}$
depend only on the residues of the numbers $q$ and $l$ modulo $8$  
We will not indicate relations of the form
$[e_i, e_j]=0$ in the definitions of Lie algebras.

The first two algebras $\mathfrak{m}_0$ and $\mathfrak{m}_2$ appeared in the work by Vergne  \cite{V}.  The Lie algebra $W^+$ is the positive part of the Witt algebra. The Lie algebras ${\mathfrak n}_1$ and ${\mathfrak n}_2$ are maximal nilpotent Lie subalgebras of the affine Kac-Moody algebras $A_1^{(1)}$ and $A_2^{(2)}$ \cite{Kac}.

\begin{table}
\caption{Structure constants of the Lie algebra ${\mathfrak n}_2$.}
\label{structure_const_n_2}
\begin{center}
\begin{tabular}{|c|c|c|c|c|c|c|c|c|}
\hline
&&&&&&&&\\[-10pt]
 &$f_{8j}$  &$f_{8j{+}1}$ &$f_{8j{+}2}$&$f_{8j{+}3}$&$f_{8j{+}4}$&$f_{8j{+}5}$&$f_{8j{+}6}$&$f_{8j{+}7}$\\
&&&&&&&&\\[-10pt]
\hline
&&&&&&&&\\[-10pt]
$f_{8i}$ & $0$ &$1$ & ${-}2$ & ${-}1$ &$0$ &$1$ &$2$ &${-}1$\\
&&&&&&&&\\[-10pt]
\hline
&&&&&&&&\\[-10pt]
$f_{8i{+}1}$ & ${-}1$ & $0$ & $1$ & $1$ & ${-}3$ & ${-}2$ & $0$ & $1$\\
&&&&&&&&\\[-10pt]
\hline
&&&&&&&&\\[-10pt]
$f_{8i{+}2}$ & $2$ & ${-}1$ &$0$ &$0$&$0$&$1$&${-}1$&$0$ \\
&&&&&&&&\\[-10pt]
\hline
&&&&&&&&\\[-10pt]
$f_{8i{+}3}$ & $1$ &${-}1$&$0$&$0$&$3$&${-}1$&$1$&${-}2$ \\
&&&&&&&&\\[-10pt]
\hline
&&&&&&&&\\[-10pt]
$f_{8i{+}4}$ & $0$ &$3$&$0$&${-}3$&$0$&$3$&$0$&${-}3$ \\
&&&&&&&&\\[-10pt]
\hline
&&&&&&&&\\[-10pt]
$f_{8i{+}5}$ & ${-}1$ &$2$&${-}1$&$1$&${-}3$&$0$&$0$&${-}1$ \\
&&&&&&&&\\[-10pt]
\hline
&&&&&&&&\\[-10pt]
$f_{8i{+}6}$ & ${-}2$ &$0$&$1$&${-}1$&$0$&$0$&$0$&$1$ \\
&&&&&&&&\\[-10pt]
\hline
&&&&&&&&\\[-10pt]
$f_{8i{+}7}$ & $1$ &${-}1$&$0$&$2$&$3$&$1$&${-}1$&$0$ \\[2pt]
\hline
\end{tabular}
\end{center}
\end{table}

\begin{theorem}[A. Fialowski, \cite{Fial}] 
Let ${\mathfrak g}{=}\oplus_{i{=}1}^{+\infty} {\mathfrak g}_i$ be a ${\mathbb N}$-graded Lie algebra of width one and let it be generated by two homogeneous components ${\mathfrak g}_1$ and ${\mathfrak g}_2$. Then ${\mathfrak g}$ is isomorphic to one and only one algebra from the following list
$$
{\mathfrak m}_0, {\mathfrak m}_2, W^+, {\mathfrak n}_1, {\mathfrak n}_2, 
\left\{{\mathfrak g}(\lambda_8, \lambda_{12},{\dots},\lambda_{4k},\dots), \lambda_{4k} \in {\mathbb P}{\mathbb K}^1, k \ge 2\right\},
$$
where $\left\{{\mathfrak g}(\lambda_8, \lambda_{12},{\dots},\lambda_{4k},\dots), \lambda_{4k} \in {\mathbb P}{\mathbb K}^1, k \ge 2\right\}$ denote the family of Lie algebras,
depending on a countable number of parameters $\lambda_{4k} \in {\mathbb P}{\mathbb K}^1, k \ge 2$.
\end{theorem}

A finite-dimensional positive-graded Lie algebra ${\mathfrak g}=\oplus_{i=1}^{+\infty}{\mathfrak g}_i$ is obviously nilpotent.

\begin{definition}
A Lie algebra $\mathfrak {g}$ is said to be pro-nilpotent if for ideals ${\mathfrak g}^i, {\mathfrak g}^1={\mathfrak g}, {\mathfrak g}^{i+1}=[{\mathfrak g}^1,{\mathfrak g}^i]$ if the ideals of its lower central series satisfy
y
$$
\cap_{i{=}1}^{\infty}\mathfrak{g}^i=\{0\}, \;\;\;\;
\dim{{\mathfrak g}/{\mathfrak g}^{i}} <+\infty,\; \forall i \in {\mathbb N}.
$$
\end{definition}

We consider the pro-nilpotent Lie algebra ${\mathfrak g}$ and its associated graded Lie algebra
${\rm gr} \mathfrak{g}=\oplus_{i=1}^{+\infty}  \left(\mathfrak{g}^i / \mathfrak{g}^{i{+}1}\right)$ with respect to the filtration by ideals of the lower central series.

\begin{definition}
A Lie algebra $\mathfrak {g}$ is said to be naturally graded if it is isomorphic to its associated graded Lie algebra
${\rm gr} \mathfrak{g}$.
\end{definition}

From the properties of the lower central series one can derive one very important property of natural grading
$
[{\mathfrak g}_1,{\mathfrak g}_i]={\mathfrak g}_{i{+}1}, i \in {\mathbb N}.
$
In particular, a naturally graded Lie algebra ${\mathfrak g}=\oplus_{i=1}^{+\infty}{\mathfrak g}_i$  is generated as a Lie algebra
by one component ${\mathfrak g}_1$.
The Lie algebra $\mathfrak {m} _0$, considered above, is naturally graded. However, its natural grading differs from its grading, which we considered in the previous section. The natural grading corresponds to such a direct sum with the two-dimensional first homogeneous subspace and shifted by one by grading
$$
{\mathfrak m}_0= \langle e_1, e_2\rangle \oplus \langle e_3 \rangle \oplus \langle e_4 \rangle\oplus\dots\oplus  \langle e_k \rangle \oplus \dots 
$$
Lie algebras ${\mathfrak n}_1$ and ${\mathfrak n}_2$ also are naturally graded
$$
\begin{array}{c}
{\mathfrak n}_1=\oplus_{i=1}^{+\infty} ({\mathfrak n}_1)_i, ({\mathfrak n}_1)_{2k+1}=\langle e_{3k+1}, e_{3k+2} \rangle,  ({\mathfrak n}_1)_{2k+2}=\langle e_{3k+3}\rangle, k \ge 0; \\
{\mathfrak n}_2=\oplus_{i=1}^{+\infty} ({\mathfrak n}_2)_i, ({\mathfrak n}_2)_{6k+1}=\langle e_{8k+1}, e_{8k+2} \rangle, ({\mathfrak n}_2)_{6k+s}=\langle e_{8k+s+1}\rangle, s=2,3,4,\\
 ({\mathfrak n}_2)_{6k+5}=\langle e_{8k+6}, e_{8k+7} \rangle, ({\mathfrak n}_1)_{6k+6}=\langle e_{8k+8}\rangle, k \ge 0.
\end{array}
$$ 
The positive part of $W^+$ of the Witt algebra and the second algebra of $\mathfrak{m}_2$ are not naturally graded Lie algebras.

The formulas
\begin{equation}
\label{real_forms}
[u,v]=w, [v,w]=\pm u, [w,u]=v
\end{equation}
determine the structure constants of two real forms of a complex simple Lie algebra ${\mathfrak sl}(2,{\mathbb C})$, the plus sign in (\ref{real_forms}) matches
the Lie algebra ${\mathfrak so}(3, {\mathbb R})$, the minus sign matches the Lie algebra ${\mathfrak so}(1,2)$. We consider two subalgebras ${\mathfrak n}_1^+$ and ${\mathfrak n}_2^-$ in loop algebras ${\mathfrak so}(3, {\mathbb R})\otimes {\mathbb R}[t]^+$ and 
${\mathfrak so}(1,2)\otimes {\mathbb R}[t]^+$
respectively. Each of these Lie subalgebras is given as the linear span of monomials of the form $\langle u \otimes t^{2k+1}, v \otimes t^{2k+1}, w \otimes t^{2k+2}, k \ge 2 \rangle$. The structure of the Lie algebra in the loop algebra
over the Lie algebra ${\mathfrak g}$ is defined in the standard way
$$
[a\otimes P(t), b \otimes Q(t)]=[a,b] \otimes P(t)Q(t), \forall a, b \in {\mathfrak g}, \forall P(t), Q(t) \in {\mathbb R}[t].
$$

\begin{lemma}
The Lie algebras ${\mathfrak n}_1^{\pm}$  are isomorphic over
${\mathbb C}$ and non-isomorphic over  ${\mathbb R}$.
\end{lemma}

\begin{definition}
Define ${\mathfrak n}_2^3$ as a one-dimensional central extension of the Lie algebra ${\mathfrak n}_2$. To set the basis and relations of ${\mathfrak n}_2^3$ we add to the infinite basis  $e_1, e_2, e_3\dots$ of the Lie algebra ${\mathfrak n}_2$ a new (central) element $z$, satisfying the additional commutation relations
$$
[e_2,e_3]=z, [z,e_i]=0,  i \in {\mathbb N}.
$$
\end{definition}
It is obvious that the Lie algebras ${\mathfrak n}_1^{\ pm}$ and ${\mathfrak n}_2^3$ are naturally graded.

For further discussion we need the apparatus of central extensions of Lie algebras generalizing an example of the Lie algebra ${\mathfrak n}_2^3$. We briefly recall the necessary definitions.

Consider the linear subspace $V$ as an abelian Lie algebra.
 The central extension of the Lie algebra
$\mathfrak{g}$ is an exact sequence
\begin{equation}
\label{exactseq}
\begin{CD}0 @>>> V@>{i}>>\tilde {\mathfrak g} @>\pi>>{\mathfrak g}@>>>0
\end{CD}
\end{equation}
of Lie algebras and their homomorphisms in which the image of the homomorphism
$i: V \to \tilde {\mathfrak{g}}$ 
is in the center $Z(\tilde{\mathfrak{g}})$ of the Lie algebra $\tilde{\mathfrak{g}}$.

As a vector space, the Lie algebra $\tilde{\mathfrak {g}}$ is a direct sum 
$V \oplus \mathfrak{g}$ with the standard inclusion $i$ and projection $\pi$.
The Lie bracket in the vector space $V\oplus \mathfrak {g}$ is defined by the formula
$$
[(v,g), (w, h)]_{\tilde {\mathfrak g}}=(c(g,h), [g,h]_{{\mathfrak g}}), \; \; \; g, h \in {\mathfrak g},
$$
where $c$ is a bilinear function on ${\mathfrak g}$, which takes its values in the space $V$.
One can verify directly that the Jacobi identity for this bracket is equivalent to the fact that the bilinear function $c$
is a cocycle, i.e. the identity holds
$$
c([g,h]_{{\mathfrak g}},e)+c([h,e]_{{\mathfrak g}},g)+c([e,g]_{{\mathfrak g}},h)=0, \;\; \forall g,h,e \in {\mathfrak g},
$$
In this connection we assume the validity of the Jacobi identity for the original Lie bracket $[g,h]_{{\mathfrak g}}$. Among the cocycles, coboundaries are distinguished $db$, skew-symmetric bilinear functions of the form $db(x,y)=b([x,y])$, where
$b$ -- some linear map ${\mathfrak g}\to V$. The quotient space $Z^2({\mathfrak g},V)/B^2({\mathfrak g},V)$ of the cocycle space over the coboundary space is called the second cohomology space $H^2({\mathfrak g},V)$ of the Lie algebra ${\mathfrak g}$ with coefficients in the trivial ${\mathfrak g}$-module $V$. 

Let $V={\mathbb K}$ be a one-dimensional space and ${\mathfrak g}=\oplus_{i=1}^{+\infty}$ be a ${\mathbb N}$-graded Lie algebra. We select in the space of cocycles $Z^2({\mathfrak g}, {\mathbb K})$ a subspace $Z^2_{(k)}({\mathfrak g},{\mathbb K})$ of cocycles with grading $k$
$$
Z^2_{(k)}({\mathfrak g},{\mathbb K})=\left\{c \in Z^2({\mathfrak g},{\mathbb K}) | c(x,y)=0, x \in  {\mathfrak g}_i,  y \in  {\mathfrak g}_j, i+j \ne k \right\}
$$
Coboundaries $B^2_{(k)}({\mathfrak g}, {\mathbb K})$ of gradings $k$ correspond to linear functions $b$ from ${\mathfrak g}^*$ that are nontrivial only on the subspace ${\mathfrak g}_k$. Corresponding subspaces $H^2_{(k)}({\mathfrak g},{\mathbb K})$ define a positive grading in the two-dimensional cohomology of the Lie algebra ${\mathfrak g}$
$$
H^2({\mathfrak g},{\mathbb K})=\oplus_{k=1}^{+\infty }H^2_{(k)}({\mathfrak g},{\mathbb K})
$$ 
The second cohomology of the Lie algebras ${\mathfrak n}_1$, ${\mathfrak n}_2$ are two-dimensional  \cite{Garl}. However the space $H^2({\mathfrak m}_0,{\mathbb K})$ is infinite-dimensional but all its homogeneous subspaces
$H^2_{(k)}({\mathfrak m}_0,{\mathbb K})$ not more than one-dimensional
$$
H^2_{(2m)}({\mathfrak m}_0,{\mathbb K})=0, \; H^2_{(2m-1)}({\mathfrak m}_0,{\mathbb K})=\langle \omega_{2m-1}=\sum_{l=2}^{m}({-}1)^l e^l{\wedge}e^{2m+1-l} \rangle,
$$
where $e^i \in {\mathfrak g}^*, i \in {\mathbb N},$  form a dual basis to the basis (\ref{5primerov}). We recall that the natural grading of the basis elements $e_i, i \ge 2$ is shifted by one unit to the left in comparison with the graduation of the width one considered above.

We consider an arbitrary subset $S$ (finite or infinite) of the set of odd positive integers greater than one. We order elements of $ S $ in ascending order $1< r_1 < r_2 < r_3 < \dots$. To every odd number
$r_j=2m_j-1$ we associate a cocycle $\omega_{r_j}=\sum_{l=2}^{m_j}({-}1)^l e^l{\wedge}e^{2m_j+1-l}$ and define a central extension (possibly infinite-dimensional) ${\mathfrak m}_0^S$ corresponding to a set of cocycles $\left\{ \omega_{r_j}, r_j \in S\right\}$.
\begin{definition}
Let $S \subset \{3,5,7,\dots\}$. Define a Lie algebra ${\mathfrak m}_0^S$
by an infinite basis consisting of the vectors $e_i, i \in {\mathbb N}$ and $z_ {r_j}, r_j \in S$ and commutation relations
\begin{equation}
[e_1,e_i]=e_{i+1}, \; i \ge 2, 
[e_l,e_{r_j-l+2}]=({-}1)^{l+1}z_{r_j}, 2\le l \le \frac{r_j+1}{2},\; r_j \in S.
\end{equation}
\end{definition}
Now everything is ready to formulate the main theorem.

\section{The formulation of the main theorem}

\begin{theorem}
\label{osnovn}
Let $\mathfrak{g}=\bigoplus_{i=1}^{+\infty}\mathfrak{g}_i$  be an infinite-dimensional ${\mathbb N}$-graded Lie algebra over ${\mathbb R}$ such that
\begin{equation}
\label{3/2}
[\mathfrak{g}_1, \mathfrak{g}_i]=\mathfrak{g}_{i+1}, \;\;
\dim{\mathfrak{g}_i}+\dim{\mathfrak{g}_{i{+}1}}\le 3, \; i \in {\mathbb N}.
\end{equation}
Then $\mathfrak{g}$ is isomorphic to one and only one Lie algebra from the following list
$$\mathfrak{n}_1^{\pm}, \mathfrak{n}_2, \mathfrak{n}_2^3, \mathfrak{m}_{0}, \{ \mathfrak{m}_{0}^{S}, S \subset \{3,5,7,\dots,2m{+}1,\dots \}.$$

The classification over ${\mathbb C}$ differs from the real classification only in that the Lie algebras $\mathfrak {n}_1^{\pm}$ are isomorphic over ${\mathbb C}$.
\end{theorem}
\begin{remark}The Lie algebras of the family ${\mathfrak g}(\lambda_{8}, \lambda_{12}, \dots) $ from Fialowski's theorem \cite{Fial}, although they are naturally graded of widths two, however in their natural grading all homogeneous components are two-dimensional for  $i \ne 2$.
\end{remark}

\section{Carnot algebras and the proof of the main theorem}
A finite-dimensional naturally graded Lie algebra ${\mathfrak g}$ is nilpotent and such a Lie algebra is said to be in the sub-Riemannian geometry by the Carnot algebra \cite{AgrM}. 
\begin{definition}
A finite-dimensional Lie algebra ${\mathfrak g}$ is called a Carnot algebra if it has ${\mathbb N}$-grading
${\mathfrak g}=\oplus_{i=1}^n{\mathfrak g}_i$ such that
\begin{equation}
\label{carnot}
[{\mathfrak g}_1,{\mathfrak g}_i]={\mathfrak g}_{i{+}1}, i=1,2,\dots,n-1, \; [{\mathfrak g}_1,{\mathfrak g}_n]=0.
\end{equation}
The maximal number $m$ for which ${\mathfrak g}_m \ne 0$ will be called the length of the Carnot algebra. It is obvious that the length of the Carnot algebra is equal to its nil-index.
\end{definition}

\begin{example}
A Lie algebra $\mathfrak{m}_0(n)$, defined by its basis
$e_1,e_2,\dots,e_n$ and commutating relations 
$[e_1,e_i]=e_{i+1}, \;2\le i \le n-1,$
is a Carnot algebra of length $n-$, this length is maximal for a given dimension $n$ and a Lie algebra with this property is called the filiform Lie algebra \cite{V} or the Lie algebra of maximal class \cite{ShZ}.
\end{example}

The proof of the main theorem will be based on three basic lemmas, which we now formulate one by one.
\begin{lemma}
Let ${\mathfrak{g}}=\oplus_{i{=}1}^{+\infty}{\mathfrak g}_i$ be an infinite-dimensional naturally graded Lie algebra. Then its arbitrary quotient Lie algebra of the form ${\mathfrak{g}}/{\mathfrak{g}}^{k+1}$
 is a Carnot algebra of length  $k$.
\end{lemma}

Let $\tilde{\mathfrak{g}}=\oplus_{i{=}1}^k{\mathfrak g}_i$ be a Carnot algebra of length $k$.
The last homogeneous summand ${\mathfrak g}_k$ belongs to the center
$Z(\tilde {\mathfrak g})$.
Consider a quotient Lie algebra ${\mathfrak g}=\tilde{\mathfrak{g}}/{\mathfrak g}_k$. It is also a Carnot algebra, which can be written in the form
direct sum of subspaces ${\mathfrak g}{=}\oplus_{i{=}1}^{k{-}1}{\mathfrak g}_i$. Its length will be one less than for the original Lie algebra.
Thus, we have a central extension of the Lie algebra ${\mathfrak g}$
$$
0 \to{\mathfrak g}_k \to  \tilde{\mathfrak g}  \to {\mathfrak g} \to 0, 
$$
which corresponds to some two-dimensional cocycle $\tilde c$ in
$H^2({\mathfrak g}, {\mathfrak g}_k)$. 
We fix a basis $e_1,\dots,e_ {j_k}$ of the subspace ${\mathfrak g}_k$. The cocycle $\tilde c$ is written in the corresponding coordinates
$
\tilde c=(\tilde c_1,\dots,\tilde c_{j_k}).
$
Each coordinate component $\tilde c_l$ of the cocycle $\tilde c$ is a
two-cocycle from the subspace
$H^2_{(k)}({\mathfrak g}, {\mathbb K})$ of grading $k$, where $l=1,\dots,j_k$.

\begin{lemma}
Let ${\mathfrak g}{=}\oplus_{i{=}1}^{k{-}1}{\mathfrak g}_i$  be a Carnot algebra of length $k-1$ and  $\tilde c_1,\dots,\tilde c_{j_k}$ a collection of $j_k$ cocycles from
$H^2_{(k)}({\mathfrak g}, {\mathbb K})$, where $j_k \le \dim{H^2_{(k)}({\mathfrak g}, {\mathbb K})}$. The Lie algebra $\tilde{\mathfrak g}{=}\oplus_{i {=}1}^{k}{\mathfrak g}_i$, defined with the help of the corresponding central extension, is a Carnot algebra of length $k$ and $\dim{ {\mathfrak g}_k}=j_k$ if and only if the cocycles $\tilde c_1,\dots,\tilde c_{j_k}$ are linearly independent in the subspace of two-dimensional cohomology $H^2_{(k)}({\mathfrak g}, {\mathbb K})\subset H^2({\mathfrak g}, {\mathbb K})$ of grading $k$.
\end{lemma}

\begin{definition}
An automorphism $f$ of the Carnot algebra ${\mathfrak g}=\oplus_{i=1}^{+\infty}{\mathfrak g}_i$ is said to be graded if it is compatible with its grading
$f({\mathfrak g}_i)={\mathfrak g}_i, \; \forall i \in {\mathbb N}$.
We denote by $Aut_{gr}({\mathfrak g})$ the subgroup of graded automorphisms in the group
$Aut({\mathfrak g})$ of all automorphisms of the Lie algebra ${\mathfrak g}$.
\end{definition}
Let ${\mathfrak g}=\oplus_{i{=}1}^n{\mathfrak g}_i$ be a Carnot algebra and  $\varphi$ be its graded automorphism. Consider its restriction
$\varphi_1=\varphi |_ {{\mathfrak g}_1}$ onto the homogeneous subspace
${\mathfrak g}_1$. Fixing some basis in ${\mathfrak g}_1 $, we can define the matrix $A$ of the mapping$\varphi_1$. Because elements of the basis ${\mathfrak g}_1$ are generators of our Carnot algebra ${\mathfrak g}=\oplus_ {i{=}1}^n{\mathfrak g}_i$, then the automorphism $\varphi$ is completely defined by the matrix $A$ and we therefore have a homomorphism
$$
\Phi: Aut_{gr}({\mathfrak g}) \to GL(q,{\mathbb K}), \; q=\dim{{\mathfrak g}_1}.
$$
The relations satisfied by the generators of ${\mathfrak g}_1$ define a system of polynomial equations for the matrix elements of an arbitrary automorphism $\varphi$, thus the group ${\rm Aut}_{gr}({\mathfrak g})$ is isomorphic to some algebraic subgroup of $GL(q,{\mathbb K})$. 
The group $Aut_{gr}({\mathfrak g})$ is always non-trivial, because it contains a one-dimensional algebraic torus ${\mathbb K}^*$ of homotheties:
$$
\phi(v)=\alpha^k v, v \in {\mathfrak g}_k, k=1,\dots,n, \alpha \in {\mathbb K}^*.
$$

Consider two central extensions $\tilde{\mathfrak g}$ and $\tilde{\mathfrak g}'$ of the same Carnot algebra ${\mathfrak g}=\oplus_{i=1}^{k{-}1} {\mathfrak g}_i$. Suppose that both extensions are Carnot algebras of length $k$.

Let there also exist an isomorphism $f: \tilde{\mathfrak g} '\to \tilde{\mathfrak g}$. Then the ideals of the lower central series satisfy the equalities $f\left(\tilde {\mathfrak g}'^i\right)=\tilde {\mathfrak g}^i$. Whence follows 
$
f\left( \tilde {\mathfrak g}_k'\right)=\tilde {\mathfrak g}_k,
$
since for the Carnot algebra the last homogeneous component coincides with the last nontrivial ideal $\tilde{\mathfrak g}_k=\tilde{\mathfrak g}^k$.
Hence, we have a commutative diagram
\begin{equation}
\label{isom}
\begin{CD}0 @>>> \tilde {\mathfrak g}_k@>{i_1}>>\tilde {\mathfrak g} @>\pi_1>>{\mathfrak g}@>>>0\\
    @AAA @AA{\Psi}A @AA{f}A  @AA{\Phi}A @AAA\\
    0 @>>> \tilde {\mathfrak g}_k'@>{i_2}>>\tilde {\mathfrak g}' @>\pi_2>>{\mathfrak g}@>>>0,
\end{CD}
\end{equation}
where we denote by the symbol $\Psi$ an isomorphism of vector spaces $\tilde {\mathfrak g}_k$ and $\tilde {\mathfrak g}_k'$. The symbol $\Phi$ denotes some automorphism of the Lie algebra ${\mathfrak g}$.
It can be shown that a graded automorphism can be considered as $\Phi$.
\begin{lemma}
Let $\{ \tilde c_1,{\dots}, \tilde c_{j_k} \}$ and $\{ \tilde c_1',{\dots}, \tilde c_{j_k}'\}$ be two sets of linear independent cocyles of  $H^2_{(k)}({\mathfrak g}, {\mathbb K})$, $j_k \le \dim{H^2_{(k)}({\mathfrak g}, {\mathbb K})}$.  
They define isomorphic central extensions of $\tilde{\mathfrak g}$ and $\tilde{\mathfrak g}'$ if and only if the linear spans $L=\langle \tilde c_1,{\dots}, \tilde c_{j_k} \rangle$ and L'=$\langle \tilde c_1',{\dots}, \tilde c_{j_k}'\rangle$ lie in one orbit of the linear action of the group of graded automorphisms ${\rm Aut}_{gr}({\mathfrak g})$ acting on the Grassmannian of $j_k$-dimensional subspaces in two-cohomology $H^2_{(k)}({\mathfrak g}, {\mathbb K})$ of grading $k$.   
\end{lemma}

In \cite{Mill2} it was shown that the problem of constructing an arbitrary $\mathbb N$-graded Lie algebra of width one can be solved by an infinite process of successive one-dimensional central extensions of some finite-dimensional graded Lie algebra of width one. Later, this idea was applied to the classification of a similar class of
$\mathbb N$-graded Lie algebras \cite{BKT}. In what way will this inductive process be applied to naturally graded Lie algebras? At each successive step, our inductive process consists of two stages. Suppose that we have a Carnot algebra ${\mathfrak g}$ of length $k-1$ constructed earlier. We study $m$-dimensional extensions leading to some Carnot algebra of length $k$.

First step: computation of the subspace $H^2_{(k)} (\mathfrak{g}, {\mathbb K})$ of the two-dimensional cohomology of the grading $k$ of the already existing Lie algebra $\mathfrak {g}$.

The second stage consists in studying the orbit space of the action of the group
${\rm Aut}_{gr}({\mathfrak g})$ of graded automorphisms on the Grassmannian of $m$-dimensional linear subspaces of the space $H^2_{(k)}(\mathfrak{g}, {\mathbb K})$.

We begin our inductive process with the three-dimensional Lie algebra ${\mathfrak m}_0(2)$. It has two different central extensions: the one-dimensional extension ${\mathfrak m}_0(3)$ and the two-dimensional extension ${\mathcal L}(2,3) \cong {\mathfrak n}(3) \cong {\mathfrak m}_0^3(3)$, where ${\mathcal L}(2,3)$ denotes a free nilpotent Lie algebra with two generators of length three, and ${\mathfrak n}_1(3)={\mathfrak n}_1/{\mathfrak n}_1^4$.

The Lie algebra ${\mathfrak m}_{0}^{2j{+}1}(2m{+}1),1 \le j \le m$
is given by a basis
$e_1, e_2,\dots, e_{2m{+}1}, z_{2j{+}1}$ and structure constants
$$
[e_1,e_i]=e_{i+1}, \quad i=2,\dots,2m,\;
[e_l,e_{2j-l+1}]=({-}1)^{l+1}z_{2j+1}, \quad l=2,\dots,j.
$$
We skip the Carnot algebras of lengths $4$ and $5$, we immediately represent a list of real Carnot algebras of length $6$
$$
{\mathfrak m}_0(6), {\mathfrak m}_0^5(6), {\mathfrak m}_0^3(6), {\mathfrak m}_0^{3,5}(6), {\mathfrak n}_{1}^{+}(6), {\mathfrak n}_{1}^{-}(6), {\mathfrak n}_{2}(6), {\mathfrak n}_{2}^3(6).
$$
This will be the "progenitors" for all infinite-dimensional naturally graded Lie algebras appearing in the formulation of the main theorem.
But the path of further proof will be rather thorny, the fact is that further along the road there will arise Carnot algebras that are not extensible to Carnot algebras with a longer length. The simplest example of such an algebra is a filiform Lie algebra ${\mathfrak m}_1(2m{+}1)$ built in \cite{V}.
\begin{example}
A Lie algebra ${\mathfrak m}_1(2m{+}1), m \ge 1,$ is given by its basis
$
e_1,{\dots},e_{2m{+}1}, z,
$
and relations
$[e_1,e_i]=e_{i+1}, 2{\le} i {\le} 2m{-}1, [e_q,e_{2m+1-q}]=({-}1)^{q+1}z, 2\le q \le m$.
\end{example}
A more complicated example of a non-extendable Carnot algebra is
\begin{example} 
A Lie algebra ${\mathfrak m}_{0,3}^{r_1,\dots,r_k}(2m{+}1)$ is given by its basis
$$
e_1,{\dots},e_{2m}, e_{2m{+}1}, e_{2m{+}2}, \; 
z_{r_1},\dots, z_{r_k}, z_{2m{-}1},
$$
where  $3 \le r_1<\dots < r_k \le 2m-3$ a set of odd positive integers, if $m \ge 3$,
and commutating relations
\begin{align*}
\label{m_0_3_S}
[e_1,e_i]=e_{i{+}1}, \; i=2, \dots, 2m{-}1, \;[e_1,z_{2m-1}]=e_{2m+1}, \;  [e_1, e_{2m+1}]=e_{2m+2},\\
[e_q,e_{2m{+}1-q}]=({-}1)^{q+1}z_{2m-1}, q=2,{\dots}, m;\\
 [e_q,e_{r_j{+}2{-}q}]=({-}1)^{q+1}z_{r_j}, \; j=1,{\dots},k, \; q=2,{\dots}, \frac{r_j{+}1}{2}; \\
[e_q,e_{2m+2-q}]=({-}1)^{q+1}(m{+}1{-}q)e_{2m+1}, \; q=2,
\dots,m;\\
[e_p,e_{2m-p+2}]=
({-}1)^{p}(p{-}2)\left(m{-}\frac{p{-}1}{2}\right)e_{2m+2},\;
p=3,\dots,m{+}1.
\end{align*}
\end{example}

Among the extensions of the Lie quotient algebras ${\mathfrak n}_1^{\pm}(n)$ and ${\mathfrak n}_2(n)$ there are also non-extendable Carnot algebras. We have as many endless chains of extensions of Carnot algebras as infinite-dimensional naturally graded Lie algebras indicated in the theorem. We, however, still exclude such extensions, which lead to two two-dimensional consecutive homogeneous components. This excludes from consideration a multiparameter Fialowski's family of Lie algebras ${\mathfrak g}(\lambda_{8k},\lambda_{12k},\dots)$  \cite{Fial}. See all the details of the proof in \cite{Mill1}.

\end{document}